\documentstyle[12pt]{amsart}
\newcommand{\norm}[1]{ \parallel #1 \parallel}
\newcommand{\beq}{\begin{equation}}
\newcommand{\eeq}{\end{equation}}
\newcommand{\bit}{\begin{itemize}}
\newcommand{\eit}{\end{itemize}}
\newtheorem{defn}{Definition}[section]
\newtheorem{prop}[defn]{Proposition}

\newtheorem{lemma}[defn]{Lemma}
\newtheorem{cor}[defn]{Corollary}
\theoremstyle{remark}
\newtheorem{rem}[defn]{Remark}

\newcommand{\bea}{\begin{eqnarray}}
\newcommand{\ena}{\end{eqnarray}}
\newcommand{\beano}{\begin{eqnarray*}}
\newcommand{\enano}{\end{eqnarray*}}
\newcommand{\bee}{\begin{enumerate}}
\newcommand{\ene}{\end{enumerate}}
\newcommand{\A}{{\cal A}}

\newcommand{\Slp}{{{\cal S}(L^p(X))}}
\newcommand{\Ao}{{\cal A}_0}

\begin{document}
\thispagestyle{empty}

\vspace*{1cm}

\begin{center}
{\Large \bf $L^p$-Spaces as Quasi *-Algebras}

\vspace{2cm}

{\large F. Bagarello }
\vspace{3mm}\\
 Dipartimento di Matematica e Applicazioni \\
Facolt\`a d' Ingegneria \\
Universit\`a di Palermo \\       \baselineskip15pt
     I-90128 - Palermo - Italy \vspace{2mm}\\

\vspace{6mm}
and

\vspace{3mm}
{\large C.Trapani}
\vspace{3mm}\\
Istituto di Fisica\\
Universit\`a di Palermo\\
\baselineskip15pt
I-90123 - Palermo - Italy\vspace{2mm}

\end{center}

\vspace*{2cm}

\begin{abstract}
\noindent
The Banach space $L^p(X,\mu)$, for $X$ a compact Hausdorff measure
space, is considered as a special kind of quasi *-algebra (called CQ*-algebra)
over the C*-algebra $C(X)$ of continuous functions on $X$. It is shown that,
for $p \geq 2$, $(L^p(X,\mu), C(X))$ is *-semisimple (in a generalized sense).
Some consequences of this fact are derived.
\end{abstract}

\vfill

\newpage
\section{Introduction}

Let $\cal A$ be a linear space and $\cal A _0$ a  $^\ast$ -algebra
contained in $\cal A$.
We say that $\cal A$ is a {\em quasi } $^\ast$-algebra with {\em
distinguished}  $^\ast$-algebra  $\cal A _0$ (or, simply, over $\cal A
_0$) if (i) the right and left  multiplications of an element of $\cal A$
and an element of $\cal A _0$ are  always defined and linear; and (ii) an
involution * (which  extends the involution of $\cal A _0$) is defined in
$\cal A$ with the property  $(AB) ^\ast =B ^\ast A ^\ast$ whenever the
multiplication is defined.

Quasi *-algebras \cite{lass1, lass2} arise in natural way as completions of
locally  convex *-algebras whose multiplication is not jointly continuous; in
this case one has to deal with topological quasi $^\ast$-algebras.

A quasi $^\ast$-algebra $(\cal A,\,\cal A_0)$ is called  {\em topological} if a
locally convex topology
$\tau$ on $\cal A$ is given such that:
\begin{itemize}
\item[(i)] the involution $A\mapsto A^\ast$ is continuous
\item[(ii)] the maps $A\mapsto AB$ and $A\mapsto BA$ are continuous  for each
$B\in \cal A_0$
\item[(iii)] $\cal A_0$ is dense in $\cal A [\tau]$.
\end{itemize}

\vspace{3mm}
In a topological quasi *-algebra the associative law holds in the following two
formulations
$$ A(BC)=(AB)C; \hspace{3mm} B(AC)=(BA)C \hspace{3mm}\forall A \in \A,
\forall B,C \in \Ao$$

\vspace{3mm}
Let $(X,\mu)$ be a measure space with $\mu$ a Borel measure on the locally
compact
Hausdorff space $X$.
As usual, we denote by $L^p(X,d\mu)$ (or simply, $L^p(X)$ if no confusion  is
possible) the Banach space of all (equivalence classes of) measurable
functions $f: \,X \longrightarrow \Bbb{C}$ such that $$ \norm{f}_p \equiv
\left( \int_X |f|^p\,d\mu\right)^{1/p} <\infty. $$ On $L^p(X)$ we consider the
natural involution $f\in L^p(X)\mapsto f^\ast \in L^p(X)$ with
$f^\ast(x)=\overline{f(x)}$.

\noindent We denote with $C_0(X)$ the C*-algebra of continuous functions
vanishing at infinity.

\noindent The pair $(L^p(X,\mu),\, C_0(X))$ provides the basic commutative
example of topological quasi $^\ast$-algebra.

{}From now on, we assume that $\mu$ is a positive measure.

\vspace{3mm}
In a previous paper \cite{bt} we introduced a particular class of topological
quasi $^\ast$-algebras, called CQ*-algebras. The definition we will give here
is not the general one, but it is exactly what we need in the commutative case
which we will consider in this paper.

A CQ*-algebra is a topological quasi $^\ast$-algebra $(\cal A,\,\cal A_0)$ with
the following properties
\begin{itemize}
\item[(i)] $\Ao$ is a C*-algebra with respect to the norm $\norm{\,}_0$ and
the involution *.
\item[(ii)]$\A$ is a Banach space with respect to the norm
$\norm{\,}$ and $\norm{A^\ast}=\norm{A}\hspace{3mm}\forall A \in A$.
\item[(iii)] $$\norm{B}_0 = \max\left\{ \sup_{\norm{A}\leq 1}\norm{AB},
\sup_{\norm{A}\leq 1}\norm{BA}\right\} \hspace{3mm}\forall B \in \Ao.$$
\end{itemize}

\vspace{3mm}
It is shown in \cite{bt} that the completion of any C*-algebra $(\Ao,
\norm{\,}_0)$ with respect to a weaker norm $\norm{\,}_1$ satisfying
\begin{itemize} \item[(i)]$\norm{A^\ast}_1=\norm{A}_1\hspace{3mm}\forall A \in
\Ao$ \item[(ii)]$\norm{AB}_1\leq\norm{A}_1\norm{B}_0\hspace{3mm} \forall A,B\in
\Ao$
\end{itemize}
is a CQ*-algebra in the sense discussed above.

\vspace{3mm}
This is the reason why both $(L^p(X,\mu),\, C_0(X))$ and
$(L^p(X,\mu),L^\infty(X,\mu))$ are CQ*-algebras.

$L^p$-spaces are examples of the $L_\rho$'s considered in \cite{zan}.
Let $\mu$ be a measure in a non-empty point set $X$ and $M^+$ be the collection
 of all
the positive $\mu-$measurable functions. Suppose that to each $f\in M^+$ it
corresponds
a number $\rho (f)\in [0,\infty]$ such that:

i)  $\rho (f)=0$ iff $f=0$ a.e. in $X$;

ii) $\rho (f_1+f_2) \leq \rho (f_1)+ \rho (f_2)$;

iii) $\rho (af)=a \rho (f) \: \: \forall a\in {\Bbb R}_+$;

iv) let $f_n \in M^+$ and $f_n \uparrow f$ a.e. in $X$. Then $\rho
(f_n)\uparrow \rho
(f)$.

Following \cite{zan}, we call $\rho$ a {\em function norm}. Let us define
$L_{\rho}$ as the set of all
$\mu-\mbox{measurable}$ functions such that $ \rho (f)<\infty $.
The space $L_{\rho}$ is a Banach space, that is it is complete,
with respect to the norm $\norm{f} \equiv \rho (|f|)$.
If the function norm $\rho$ satisfies the additional condition
$$ \rho (|fg|) \leq \rho (|f|) \norm{g}_\infty, \quad \forall f,g \in
C_0(X),$$ then the completion of $C_0(X)$ with respect to this norm is an
abelian CQ*-algebra over $C_0(X)$.

Of course, for $L^p$-spaces,  $\norm{\, } = \norm{\,
}_p$.

\vspace{3mm}
In this paper we will discuss
some structure properties of $(L^p(X,\mu),\, C_0(X))$ as a CQ*-algebra. In
Section 2, in particular, we will study a certain class of positive
sesquilinear forms on
$(L^p(X,\mu),\, C_0(X))$ which lead, in rather natural way, to a definition of
{\em *-semisimplicity}. As is shown in \cite{bt2}, *-semisimple CQ*-algebras
behave
nicely and for them a refinement of the algebraic structure of quasi
$^\ast$-algebra to a {\em partial *-algebra} \cite{ait1,ait2} is possible.
The abelian case is discussed in Section 3.

Finally, we characterize *-semisimple abelian CQ*-algebras as a CQ*-algebra
of functions obtained by means of a family of $L^2$-spaces, generalizing in
this way the concept of Gel'fand transform for C*-algebras.

\section{Structure properties of $L^p$-spaces}
\begin{lemma}
Let $g$ be a measurable function on $X$, with $\mu (X)< \infty$, and  assume
that $fg \in L^r(X)$ for  all $f \in L^p(X)$ with $1\leq r\leq p$. Then $g \in
L^q(X)$ with $ p^{-1}+q^{-1}= r^{-1}$.
\label{PRD}
\end{lemma}
\begin{pf}
Let us consider the linear operator $T_g$ defined in the following way:
$$ T_g: f\in L^p(X) \longrightarrow fg \in L^r(X).$$
$T_g$ is closed. Indeed, let $f_n \stackrel{p}\longrightarrow f$ and $T_g f_n
\stackrel{r}\longrightarrow h$; this implies the existence of a subsequence
$f_{n_k}$ such that $f_{n_k} \longrightarrow f$ a.e.
in $ X$  and so $f_{n_k}g \longrightarrow fg$ a.e. in $ X$. Then, we have
$f_{n_k}g \longrightarrow fg$ in measure; on the other hand the fact that
$T_g f_n\stackrel{r}\longrightarrow h$ implies also that $f_n g$ converges to
$h$
in measure. Thus, necessarily, $h(x)=f(x)g(x)$ a.e. in $X$.

\noindent Thus $T_g$ is closed and everywhere defined in $L^p(X)$; then, by
the closed graph theorem, $T_g$ is bounded; i.e., there exists $C>0$
(depending on $g$) such that
$$ \norm{fg}_r \leq C \norm{f}_p $$
If $r=1$, this already implies that $g \in L^{p'}(X)$ with $p^{-1}+p^{'-1}=1$.

\noindent Let now $r>1$ and let $h\in L^{r'}(X)$ with $r^{-1}+r^{'-1}=1$
If $f\in L^p(X)$ then, by H\"older inequality, $fh \in L^m(X)$ with $ m^{-
1}=p^{-1}+r^{'-1}$. Thus, the set
$$ \cal F = \{ fh\, :\; f\in L^p(X), \;h\in L^{r'}(X) \}$$
is a subset of $L^m(X)$.
Conversely, any function $\psi \in L^m(X)$ can be factorized as the product of
a function in $L^p(X)$ and one in $L^{r'}(X)$. This is achieved by
considering, for instance,
the principal branches of the functions $\psi ^{m/p}$ and $\psi ^{m/{r'}}$.

\noindent We can now apply the first part of the proof to conclude that $g \in
L^{m'}(X)$ with $m^{-1}+m^{'-1}=1$. An easy computation shows now that $m'
=q$.

\end{pf}
\begin{rem}
For $r=1$ and $\mu (X)$ not necessarily finite, the above statement can also
 be found in \cite{zan,HLP}.
Using this fact, with the same technique as above, the statement of Lemma
\ref{PRD} can be extended to the case $\mu (X) =\infty$.
Our proof involves functional aspects of $L^p$-spaces
so we think is worth giving it.
\end{rem}
 \begin{defn}Let $(\A,\Ao)$ be a CQ*-algebra.
We denote as $\cal S(\A) $ the set of sesquilinear forms $\Omega$ on $\A \times
\A$ with the following properties
\begin{itemize}
\item[(i)]$ \Omega(A,A) \geq 0 \hspace{3mm} \forall A \in \A$
\item[(ii)] $\Omega(AB,C) = \Omega (B,A^*C)  \hspace{3mm} \forall A \in
\A,\hspace{2mm} \forall B,C \in \Ao$
\item[(iii)]$ |\Omega(A,B)| \leq \norm{A}
\norm{B} \hspace{3mm} \forall A,B\in \A$ \end{itemize}
The CQ*-algebra $(\A,\Ao)$ is called \underline{*-semisimple}  if
$\Omega(A,A)=0 \hspace{2mm} \forall \, \Omega \in \cal S (\A)$ implies $A=0$.
\label{SLP}
\end{defn}
{}From now on, we will only consider the case of compact $X$ (in this case we
denote as $C(X)$ the C*-algebra of continuous functions on $X$)
and
we will
focus our attention on the question whether $(L^p(X),C(X))$,
is or
is not *-semisimple. To this aim, we need first to describe the set $\Slp$.

In what follows we will often use the following fact, which we state as a lemma
for reader's convenience.
\begin{lemma}
Let $p,q,r \geq 1$ be such that $p^{-1}+q^{-1}=r^{-1}$. Let $w \in L^q(X)$.
Then the linear operator $T_w : f\in L^p(X) \mapsto fw \in L^r(X)$ is bounded
and $\norm{T_w}_{p,r}= \norm{w}_q$.
\label{normop}
\end{lemma}
\begin{rem}
 Here $\norm{T_w}_{p,r}$ denotes the norm
of $T_w$ as bounded operator from $L^p(X)$ into $L^r(X)$.
\end{rem}

For shortness, we put (if $p=2$, we set $\frac{p}{(p-2)}=\infty$).
$$ {\cal B}^p_+ = \left\{\chi  \in L^{p/(p-2)},\: \: \chi\geq
0 \mbox { and }\norm{\chi}_{p/(p-2)}\leq 1 \right\}.$$
\begin{prop} \begin{enumerate} \item
If $p\geq 2$, then any $\Omega \in \Slp$ can be represented as
\begin{equation}
 \Omega (f,g) = \int_X f(x) \overline{g(x)} \psi(x)\,d\mu
\label{sta}
\end{equation}
for some $\psi \in {\cal B}^p_+$.

Conversely,if $\psi \in {\cal B}^p_+$, then the sesquilinear form  $\Omega$
defined
by Eqn. \eqref{sta} is in $\Slp$. \item If $1\leq p<2$ then $\Slp =\{0\}$
\end{enumerate}
\label{rep}
\end{prop}
\begin{pf}
 We notice, first, that any bounded sesquilinear form $\Omega$ on $L^p(X)
\times L^p(X)$ can be represented as
\begin{equation} \Omega(f,g) = <f,Tg> =\int_X f(x) \overline{(Tg)(x)}\,d\mu
\label{SQF}
\end{equation}
 where $T$ is a
bounded linear operator from $L^p(X)$ into its dual space $L^{p'}(X),
\hspace{3mm}  p^{-1} + p^{'-1}=1$ \cite[\S 40]{Kot}.
{}From (ii) of Definition \ref{SLP} and from Eqn. \eqref{SQF} it follows easily
that $$ Tg = gTu \hspace{3mm} \forall g\in L^p(X)$$ where $u(x)=1 \;\,\forall
x\in X$.
Set $Tu=\psi$; from (i) of Definition \ref{SLP}, we get $\psi \geq 0$.

\vspace{3mm}
\noindent(1) If $p \geq 2$ then by Lemma \ref{PRD}, we get $\psi \in
L^{p/(p-2)(X)}$. Making use of Lemma \ref{normop} it is also easy to check that
$\norm{\psi}_{p/(p-2)}\leq 1.$

\vspace{3mm}
\noindent(2) Let now $1\leq p <2$ and $\psi \neq 0$. Since $\psi \geq 0$, we
can choose $\alpha >0$ in such a way that the set $Y= \{x\in X
:\,\psi (x)>\alpha\}$ has positive measure. Let $f \in L^p(Y)\setminus L^2(Y)$
(such a function always exists because of the assumption on $p$). Now define
$$\tilde{f}(x)=\left\{  \begin{array}{ll}f(x)&\mbox{if $x\in Y$}\\ 0 &
\mbox{if $x \in X\setminus Y$} \end{array}  \right. $$  Clearly, $\tilde{f}
\in L^p(X)$.

\noindent Now,
$$\Omega(f,f)=\int_X |\tilde{f}(x)|^2\psi(x)d\mu (x)=\int_Y |f(x)|^2\psi(x)d\mu
(x)  \geq \alpha \int_Y |f(x)|^2d\mu (x)= \infty$$
and this is a contradiction.
\end{pf}
\begin{rem}
{}From the above representation theorem, it follows easily that if $\Omega \in
\Slp$, then the sesquilinear form $\Omega ^*$ defined by $\Omega
^*(f,g)=\Omega(g^*,f^*)$ also belongs to $\Slp$.
\end{rem}
\begin{prop}Let $p\geq 2$.  $(L^p(X,d\mu), C(X))$ is *-semisimple.
\label{2.3}
\end{prop}
\begin{pf}
We show first that $\forall f\in L^p(X)$, there exists $\tilde\Omega \in \Slp $
such that $\tilde\Omega(f,f)=\norm{f}_p^2$.

\noindent This is achieved by setting
$$\tilde\Omega (g,h) = \norm{f}_p^{2-p}\int_X g \bar{h}|f|^{p-2}\,d\mu$$
$\tilde\Omega \in \Slp$ since the function $\psi = |f|^{p-2}\norm{f}_p^{2-p}$
belongs to the set ${\cal B}^p_+$. A direct calculation shows
that $\tilde\Omega(f,f)=\norm{f}_p^2$.

 Let us now suppose that $\Omega (f,f)=0 \: \: \forall \Omega \in
\Slp$. Then, in particular, $\tilde \Omega (f,f)=\norm{f}_p^2=0$.
Therefore $f(x)=0$ and
so $(L^p(X,\mu), C(X))$ is *-semisimple.
\end{pf}

Positive sesquilinear forms which are normalized (in the sense that
$\Omega(\Bbb{I}, \Bbb{I}) =\norm{\Bbb{I}}^2$) could be expected to play,
in this framework, the same role as {\em states} on a C*-algebra.
Indeed, for such an $\Omega$ we have
$$ \norm{\Omega}^2 \equiv \sup_{A\in\A}\frac{\Omega(A,A)}{\norm{A}^2} =1$$
The next Proposition shows, however, an essential difference between the two
frameworks.

\begin{prop}
In the CQ*-algebra $(L^p(X,\mu), C(X))$, $p\geq 2$,
there exists one and only one $\Omega \in {\cal S} (L^p(X))$ such that
$\Omega(u,u) =\norm{u}_p^2$, where $u(x)=1 \,\forall x\in X$.
\label{2.10}
\end{prop}
\begin{pf}
The sesquilinear form
$$\Omega _0(f,g)\equiv \mu (X)^{(2-p)/p}\int_X f(x) \overline{g(x)} \,d\mu$$
belongs to $\Slp$ and obviously satisfies the condition $\Omega _0 (u,u) =
\norm{u}^2_p$.

\noindent It remains to prove its uniqueness.

\noindent Let $\Omega$ satisfy the assumptions of the Proposition.

 By Proposition \ref{rep}, there exists
$\psi \in {\cal B}^p_+$ such that
$$
 \Omega (f,g) = \int_X f(x) \overline{g(x)} \psi(x)\,d\mu.
$$
Therefore,
$$ \Omega(u,u)= \int_X \psi (x)dx = \norm{\psi}_1 $$
and
$$ \norm{u}_p^2 = \mu (X)^{2/p}$$
and so we must have
\begin{equation}
\norm{\psi}_1 =\mu (X)^{2/p}.
\label{N1}
\end{equation}
On the other hand, since $\psi \in  {\cal B}^p_+$, using the inequality
$
\norm{\psi}_1 \leq \mu (X)^{2/p}\norm{\psi}_{p/(p-2)},
$ we conclude also that
\begin{equation}\norm{\psi}_{p/(p-2)} =1.
\label{N2}
\end{equation}
  We will
now prove that there exists only one $\psi$ satisfying
both \eqref{N1} and \eqref{N2}.
To show this, let us define on $X$ a new measure $\nu$ by $\nu (Y) :=
\frac{\mu (Y)}{\mu (X)}$ for any $\mu$-measurable set $Y\subseteq X$; then $\nu
(X) =1$
and $ \norm{\psi}_{1,\nu} =\norm{\psi}_{p/(p-2),\nu} = \mu (X)^{(2-p)/p}$
where $\norm{\,}_{r,\nu}$ denotes the $L^r$-norm with respect to the measure
$\nu$. From this it follows \cite{HLP} that $\psi$ is constant $\nu$-a.e. It is
then easy
to prove that $\psi (x) = \mu (X)^{(2-p)/p}$ $\mu$-a.e. in $X$.
\end{pf}

We notice that a statement analogous to Proposition \ref{2.10} does not hold
for C*-algebras where the set of states is, in general, quite rich.

\begin{rem}
If a CQ*-algebra $(\A,\Ao)$ has only one normalized positive sesquilinear
form $\Omega _0$, then it is possible to identify easily
$\A/ \mbox{Ker}\, \Omega _0$ with a linear subspace of the dual $\A'$ of $\A$
by
$$
[A] \in \A/ \mbox{Ker}\, \Omega _0 \longrightarrow F_A \in \A'
$$
where $F_A(B) \equiv \Omega_0 (B,A) \hspace{3mm} \forall B \in \A$.
The known imbeddings of $L^p$-spaces on sets of finite measures, for $p \geq
2$, provides examples of this situation (in this case $\mbox{Ker}\, \Omega
_0=\{0\}).$
\end{rem}
In \cite{bt2} we have introduced some norms on a semisimple CQ*-algebra which
play an interesting role. Their definitions in the case of $(L^p(X,\mu), C(X))$
reads as follows

\begin{equation}
\norm{f}_\alpha = \sup\{\Omega(f,f), \Omega \in \Slp\,\}
\end{equation}
and
\begin{equation}
\norm{f}_\beta = \sup\{|\Omega(f\phi,\phi)|,\hspace{1mm} \Omega \in
\Slp\;,\phi\in C(X),\, \norm{\phi}_\infty\leq 1\}
\end{equation}

\begin{prop}Let $f\in L^p(X)$, $p\geq 2$. Then
\begin{equation}
\norm{f}_\alpha = \norm{f}_p \hspace{4mm}\norm{f}_\beta \leq \norm{f}_{p/2}
\end{equation}
If $f\geq 0$, then $\norm{f}_\beta = \norm{f}_{p/2}$.
\end{prop}
\begin{pf}
To prove the equality $\norm{f}_\alpha = \norm{f}_p$,
it is enough recalling that, as shown in the proof of Proposition \ref{2.3},
there exists a  sesquilinear form $\tilde\Omega$ such that $\tilde\Omega(f,f)=
\norm{f}_p^2$. Therefore $\norm{f}_\alpha^2\geq
\tilde\Omega(f,f)=\norm{f}_p^2$. The converse inequality follows from (iii) of
Definition \ref{SLP}.

\noindent By Proposition \ref{rep} we get
\begin{eqnarray*}
\norm{f}_\beta &=& \sup\{|\Omega(f\phi,\phi)|,\hspace{1mm} \Omega \in
\Slp\;,\phi\in C(X),\, \norm{\phi}_\infty\leq 1\} \\ &=&
\sup\{\left|\int_Xf|\phi|^2\psi\,d\mu\right| ,\hspace{1mm} \psi \in {\cal
B}^p_+
\;,\phi\in C(X),\, \norm{\phi}_\infty\leq 1\} \leq \norm{f}_{p/2}
\end{eqnarray*}
On the other hand, for $f\geq 0$ we get
\begin{eqnarray*}
\norm{f}_{p/2}&=& \sup\{\left|\int_Xf \overline{\psi}\,d\mu \right|
,\hspace{1mm}\norm{\psi}_{p/(p-2)}\leq 1\} \\
&=&\sup\{\left|\int_Xf|\phi|^2\psi\,d\mu\right| ,\hspace{1mm}
\norm{\psi}_{p/(p-2)}\leq
1 \;,\phi\in C(X),\, \norm{\phi}_\infty\leq 1\}\\
&=&\sup\{\int_X f|\phi|^2|\psi|\,d\mu| ,\hspace{1mm}
\norm{\psi}_{p/(p-2)}\leq
1 \;,\phi\in C(X),\, \norm{\phi}_\infty\leq 1\}\\
&=&\sup\{|\Omega(f\phi,\phi)|,\hspace{1mm} \Omega \in
\Slp\;,\phi\in C(X),\, \norm{\phi}_\infty\leq 1\}=\norm{f}_\beta
\end{eqnarray*}
\end{pf}

Apart from $\Ao$, a *-semisimple CQ*-algebra has another distinguished subset,
denoted as $\A_\gamma$ which play an interesting role for what concerns the
functional calculus and representations in Hilbert space.
We give here the definition in the case $\A = L^p(X)$ with $p\geq 2$.
The general definition is an obvious extension of this one.
We denote as $\left(L^p(X)\right)_\gamma$ the set of all $f\in L^p(X)$
such that
\begin{equation}
\norm{f}_\gamma ^2:=\sup\left\{\frac{\Omega(f\phi,f\phi)}{\Omega(\phi,\phi)},
\hspace{2mm} \Omega \in \Slp, \phi \in C(X),\; \Omega(\phi,\phi)\neq 0
\right\}\,
< \infty
\end{equation}
Since $\forall \phi (x) \in C(X)$  and
$\forall \Omega \in \Slp$ the sesquilinear form $\Omega_{\phi}$, defined by $
\Omega_{\phi}(f,g) \equiv \displaystyle{\frac{\Omega (f\phi ,g\phi
)}{\norm{\phi}_\infty^2}}\: \forall f,g \in L^p(X)$,
also belongs to $\Slp$,
it turns out that
\begin{equation} \norm{f}_\gamma
^2=\sup\left\{\frac{\Omega(f,f)}{\Omega(u,u)},\hspace{2mm}  \Omega \in
\Slp\right\},  \end{equation}
where $u(x)=1 \;\forall x\in X$.

\begin{prop}
$\left(L^p(X)\right)_\gamma = L^\infty(X)$ and $\norm{f}_\gamma = \norm{f}_
\infty$
\label{gam}
\end{prop}
\begin{pf}
{}From the previous discussion and from Proposition \ref{rep}, it follows that
$f
\in  \left(L^p(X)\right)_\gamma$ if, and only if, $$ \norm{f}_\gamma ^2 =
\sup_{\psi \in {\cal B}^p_+} \frac{\int_{X}|f|^2\psi
d\mu}{\norm{\psi}_1} < \infty $$
This means that there exists a constant $C>0$ such that
$$
\omega_f(\psi) := \int_{X}|f|^2\psi \,d\mu \leq C \norm{\psi}_1\hspace{3mm}
\forall \psi \in  {\cal B}^p_+$$
Let now $\psi \in L^{p/(p-2)}(X), \: \psi \geq 0$. Then the function
$\psi_N \equiv \displaystyle{\frac{\psi}{\norm{\psi}_{p/(p-2)}}} \in {\cal
B}^p_+$. This implies that the
above inequality also holds for all positive functions in $L^{p/(p-2)}(X)$.
Finally, if $\psi$ is an arbitrary function in $L^{p/(p-2)}(X)$, then
\begin{eqnarray*} |\omega_f(\psi)|&=&\left| \int_{X} |f|^2\psi d\mu \right|
\leq \int_{X} |f|^2|\psi| d\mu \\&\leq& C \norm{|\psi|}_1= C \norm{\psi}_1.
\end{eqnarray*}
Therefore $\omega_f$
is a linear functional on $L^{p/(p-2)}(X)$ and it is bounded with respect to
$\norm{\;}_1$. Since $L^{p/(p-2)}(X)$ is dense in $L^1(X)$, $\omega_f$ has a
unique continuous extension to $L^1(X)$. Then $|f|^2 \in \L^\infty(X)$; this,
in turn, implies that $f\in L^\infty(X)$. Thus $\left(L^p(X)\right)_\gamma
\subseteq L^\infty (X)$.
\noindent
 The converse inclusion is obvious.
\noindent
 The norm of $\omega_f$ as linear functional on the subspace
$L^{p/(p-2)}(X)$ of $L^1(X)$ is
$\norm{\omega_f} = \norm{|f|^2}_{\infty} = \norm{f}^2_{\infty}$ (the latter
equality follows from the C*-nature of $L^{\infty}(X))$. Moreover, the Hahn-
Banach theorem, used to extend $\omega_f$ from $L^{p/(p-2)}(X)$ to $L^1(X)$,
ensures that $\norm{f}^2_{\gamma} = \norm{\omega_f}$.
\end{pf}

The role of $\left(L^p(X)\right)_\gamma$ will be clearer if we consider the
GNS-construction of an abstract CQ*-algebra $(\A, \Ao)$ obtained via a
sesquilinear form $\Omega$ in ${\cal S}(\A)$. This problem was, from a general
point of view, considered in \cite{bt}. We will give here a simplified version
of the GNS-construction which is closer to that proved in \cite{ait2} for
general partial *-algebras.

Let $(\A, \Ao)$ be a CQ*-algebra and $\Omega$ a positive sesquilinear form in
${\cal S}(\A)$. Let ${\cal K}= \{A \in \A:\, \Omega (A,A)=0\}$. Let us consider
the linear space $\A / {\cal K}$; an element of this set will be denoted as
$\lambda_{\Omega}(A),\;A\in \A$. Clearly, $\A / {\cal K}=\lambda_{\Omega}(\A)$
is a pre-Hilbert space with respect to the scalar product
$(\lambda_{\Omega}(A),\lambda_{\Omega}(B))= \Omega (A,B), \; A,B \in \A$. We
denote by ${\cal H}_{\Omega}$ the Hilbert space obtained by the `
completion of $\lambda_{\Omega}(\A)$. Then $\Omega$ is {\em invariant} in the
sense of \cite{ait2}. This means, in this case, that $\Omega$ satisfies
condition (ii) of Definition \ref{SLP} and that $\lambda_{\Omega}(\Ao)$ is
dense in ${\cal H}_{\Omega}$. Indeed, let $\lambda_{\Omega}(A) \in
\lambda_{\Omega}(\A)$ and let $\{A_n\}$ be a sequence in $\Ao$ converging to
$A$ in the norm of $\A$. Then from
$$\Omega(A-A_n, A-A_n) \leq \norm{A-A_n}^2$$ it follows that
$\lambda_{\Omega}(A_n)\rightarrow \lambda_{\Omega}(A)$ in ${\cal H}_{\Omega}$.

\noindent If we put
$$ \pi_{\Omega}(A)\lambda_{\Omega}(B)= \lambda_{\Omega}(AB)\hspace{3mm} B\in
\Ao,$$
then $\pi_{\Omega}(A)$ is a well-defined closable operator with domain
$\lambda_{\Omega}(\Ao)$ in ${\cal H}_{\Omega}$. More precisely it is an element
of the partial O*-algebra ${\cal L}^+(\lambda_{\Omega}(\Ao),{\cal H}_{\Omega})$
\cite{ait1, ait2}. The map $A \mapsto \pi _{\Omega}(A)$ is a *-representation
of partial *-algebras in the sense of \cite{ait2}.

\noindent We define now the following set:
$$ {\cal D}_{\Omega}=\left\{ A\in \A :\; \sup_{B\in \Ao, \Omega (B,B) \neq 0}
\frac{\Omega (AB,AB)}{\Omega (B,B)} < \infty\right\}$$
then
\begin{itemize}
\item[(i)] ${\cal D}_{\Omega}$ is a linear space;
\item[(ii)] ${\cal D}_{\Omega} \supset \Ao$;
\item[(iii)] if $A\in {\cal D}_{\Omega}$ and $ B\in \Ao$, then $AB \in {\cal
D}_{\Omega}$
\end{itemize}
If ${\cal D}_{\Omega}= \A$ then $\Omega$ is {\em admissible} in the sense of
\cite{bt}.

{}From the definition itself, it follows easily that $ \pi_{\Omega}({\cal
D}_{\Omega}) \subseteq {\cal B}({\cal H}_{\Omega})$, i.e. each element of
${\cal D}_{\Omega}$ is represented by a bounded operator in Hilbert space.

\noindent As an example, let us consider the CQ*-algebra $(L^p(X,\mu), C(X))$
for $p\geq 2$. Let us fix $w \in {\cal B}_+^p$ with $w>0$ and define
$$ \Omega (f,g)= \int_X f \overline g w d\mu \hspace{3mm} f,g \in L^p(X).$$
Thus $\Omega \in \Slp$.
It is clear that ${\cal H}_{\Omega} = L^2(X, \mu _w)$ where $d\mu_w =wd\mu$.
The representation $\pi _{\Omega}$ is then defined by
$$ \pi _{\Omega}(f)\phi = f\phi \hspace{3mm} \phi \in C(X)$$
for $f \in L^p(X)$.
A proof analogous to that of Proposition \ref{gam} shows that, in this case,
${\cal D}_{\Omega}= L^{\infty}(X,\mu_w) \cap L^p(X,\mu)$.

\begin{rem}
The set $\A_\gamma$ is included in ${\cal D}_{\Omega} \: \forall \Omega \in
{\cal S}(\A)$. This makes clear the role of $\A_\gamma$ as a relevant subset of
$\A$ represented by bounded operators $\forall \Omega \in {\cal S}(\A)$.
\end{rem}

\vspace{3mm}
We conclude this Section by showing that each *-semisimple abelian CQ*-
algebra, can be thought as CQ*-algebras of functions.

Let $X$ be a compact Hausdorff space and ${\cal M}= \{ \mu_\alpha; \, \alpha
\in {\cal I} \}$ be a family of Borel measures on $X$. Let us assume that there
 exists $C>0$ such that $\mu_\alpha (X) \leq C \; \forall \alpha
\in {\cal I}$. Let us denote by $\norm{\,}_{p,\alpha}$ the norm in $L^p(X,
\mu_\alpha)$ and define, for $\phi \in C(X)$
$$ \norm{\phi}_{p,{\cal I}} \equiv \sup_{\alpha \in {\cal I}}
\norm{\phi}_{p,\alpha}.$$
Since $\norm{\phi}_{p,{\cal I}} \leq C \norm{\phi}_\infty \hspace{3mm} \forall
\phi \in C(X)$, then $\norm{\,}_{p,{\cal I}}$ is finite on $C(X)$ and really
defines a
norm on $C(X)$ satisfying
\begin{itemize}
\item[(i)]$\norm{\phi^*}_{p,{\cal I}}=\norm{\phi}_{p,{\cal I}}\hspace{3mm}
\forall \phi \in C(X)$
\item[(ii)]$\norm{\phi \psi}_{p,{\cal I}}\leq \norm{\phi}_{p,{\cal I}}
\norm{\psi}_\infty \hspace{3mm} \forall \phi,\psi \in C(X).$
\end{itemize}
Therefore the completion $L^p_{\cal I} (X, {\cal M})$ of $C(X)$ with respect
to $\norm{\,}_{p,{\cal I}}$ is an abelian CQ*-algebra over $C(X)$.
It is clear that $L^p_{\cal I} (X, {\cal M})$ can be identified with a
subspace of $L^p(X, \mu_\alpha) \;\forall\alpha \in {\cal I}$.

\noindent It is obvious that $L^p_{\cal I} (X, {\cal M})$ contains also non-
continuous functions.

For $p\geq 2$, the CQ*-algebra $\left( L^p_{\cal I} (X, {\cal M}), C(X)
\right)$ is *-semisimple (this depends on the fact that each element of
$\Slp$ gives rise, by restriction, to an element of ${\cal S}\left( L^p_{\cal
I}
(X, {\cal M})\right)$).

\begin{prop}
Let $(\A,\Ao)$ be a *-semisimple abelian CQ*-algebra with unit $\Bbb{I}$.
Then there exists a
family ${\cal M }$ of Borel measures on the compact space X of characters of
$\Ao$ and a map
$\Phi: A\in \A \rightarrow
\Phi (A)\equiv \hat{A} \in L^2_{\cal I} (X, {\cal M})$ with the properties
\begin{itemize}
\item[(i)]
$\Phi $ extends the Gel'fand transform of elements of $\Ao$ and $\Phi
(\A) \supseteq C(X)$
\item[(ii)]
$\Phi $ is linear and one-to-one
\item[(iii)]
$\Phi  (AB) = \Phi  (A) \Phi  (B)  \hspace{3mm} \forall
A\in \A, \,B \in \Ao.$
\item[(iv)]
$\Phi  (A^*)=\Phi (A) ^* \hspace{3mm} \forall A\in \A.$
\end{itemize}
Thus $\A$ can be identified with a subspace of $L^2_{\cal I} (X, {\cal M})$

If $\A$ is \underline{regular}, i.e. if
$$ \norm{A}^2 =\sup_{\Omega \in {\cal S}(\A)} \Omega (A,A)$$
then $\Phi$ is an isometric *-isomorphism of $\A$ onto $L^2_{\cal I} (X, {\cal
M})$.
\end{prop}

\begin{pf}
Define first $\Phi $ on $\Ao$ as the usual Gel'fand transform
$$\Phi:\, B\in \Ao \rightarrow \hat{B} \in C(X)$$
where $X$ is the space of characters of $\Ao$.

As is known, the Gel'fand transform  is an isometric *-isomorphism of $\Ao$
onto $C(X)$.

\noindent Let $\Omega \in {\cal S}(\A)$ and define the linear
functional $\omega$ on $C(X)$ by
$$
\omega (\hat{B}) = \Omega  (B,\Bbb{I}).
$$
It is easy to check that $\omega$ is bounded on $C(X)$; then by the Riesz
representation theorem, there exists a unique positive Borel measure $\mu
_{\Omega}$ on $X$ such that
$$
\omega (\hat{B}) = \Omega  (B,\Bbb{I})= \int_X \hat{B} (\eta)d \mu _{\Omega}
(\eta) \hspace{3mm} \forall B \in \Ao.
$$
We have $\mu _{\Omega}(X) \leq \norm{\Bbb{I}}^2\;\, \forall \Omega \in
{\cal S}(\A)$

\noindent Let ${\cal M} \equiv \{ \mu_\Omega :\; \Omega \in {\cal S}(\A) \}$
and let $L^2_{{\cal S}(\A)}(X, {\cal M})$ be the CQ*-algebra constructed as
above.
Now, if $A \in \A$ there exists a sequence $\{ A_n\} \subset \Ao$ converging to
$A$ in the norm of $\A$. We have then
$$
\norm{ \hat{A_n} - \hat{A_m}}^2_{2,{\cal S}(\A)}
= \sup_{\Omega \in {\cal S}(\A)} \Omega(A_n - A_m, A_n - A_m) \leq \norm{A_n -
A_m}^2 \rightarrow 0
$$
Let $\hat{A}$ be the $\norm{\,}_{2,{\cal S}(\A)}$-limit in  $L^2_{{\cal
S}(\A)}(X, {\cal M})$
of $\left\{ \hat{A_n}\right\}$ and define
$$ \Phi (A) = \hat{A}.$$
Evidently, $\norm{\hat{A}}_{2,{\cal S}(\A)} = \sup_{\Omega \in {\cal S}(\A)}
\Omega (A,A)$.
This implies that if $\hat{A}=0$, then $\Omega (A,A)=0\;\,\forall \Omega \in
{\cal S}(\A)$ and thus $A=0$ for $\A$ is *-semisimple.
The proof of (ii), (iii) and (iv) is straightforward.

Now if $\A$ is regular, from the above discussion it follows immediately that
$\Phi$ is an isometry.
We conclude by proving that in this case $\Phi$ is onto. Let $\hat{A}$ be an
element of
$L^2_{{\cal S}(\A)}(X, {\cal M})$. Then there exists a sequence $\{ \hat{A}_n\}
\subset C(X)$
converging to $\hat{A}$ with respect to $\norm{\,}_{2,{\cal S}(\A)}$. Then the
corresponding sequence $\{ A_n\} \subset \Ao$ converges to $A$ in the norm of
$\A$, since $\Phi$ is an isometry.
\end{pf}

\section{The partial multiplication}
In this final Section, we will discuss the possibility of refining the
multiplication structure of $L^p$-spaces. Actually, we will show that $L^p(X)$
is really a {\em partial *-algebra}. For reader's convenience we repeat here
the definition.

\vspace{3mm}
A partial *-algebra is a vector space $\cal A$ with involution $A \rightarrow A
^\ast$
[i.e. $(A+\lambda B) ^\ast =A ^\ast +\overline{\lambda} B ^\ast$ ; $A=A^{\ast
\ast}$ ] and a
subset  $\Gamma \subset\cal A\times\cal A$ such that (i)  $(A,B)\in  \Gamma$
implies $(B ^\ast ,A
^\ast )\in  \Gamma$ ; (ii) $(A,B)$ and $(A,C)\in  \Gamma$  imply  $(A,B+\lambda
C)\in  \Gamma$ ;
and (iii) if $(A,B)\in  \Gamma$ , then there exists an element  $AB\in \cal A$
and for this
multiplication the distributive property holds  in the following sense: if
$(A,B)\in  \Gamma$  and
$(A,C)\in  \Gamma$  then

$$AB+AC=A(B+C)$$

Furthermore  $(AB) ^\ast =B ^\ast A ^\ast$ .

The product is not required to be associative.

The partial $^{*}-$algebra $\cal A$ is said to have a unit if there
exists an element $\Bbb{I}$ (necessarily unique) such that $\Bbb{I}^\ast
=\Bbb{I},\;
(\Bbb{I} ,A)\in  \Gamma,\;  \Bbb{I} A=A\Bbb{I},\; \forall A\in \cal A$.

If $(A,B)\in  \Gamma$  then we say that $A$ is a left multiplier of $B$ [and
write $A\in L(B)$] or $B$ is a right multiplier of $A$ [$B\in R(A)$].

There are at least two ways to introduce partial multiplications in $L^p(X)$.
The first one is almost obvious, if the set of multicable elements is defined
as follows:

$$\Gamma_1 =\{(f,g) \in L^p(X) \times L^p(X)\,:\, fg \in L^p(X) \}.$$

The second \cite{ait1} is obtained by defining the following set of real
numbers
$$ E(f) = \{q\in [1,\infty)\,:\, \norm{f}_q<\infty \}.$$

The partial multiplication is then defined on the set
$$ \Gamma_2=\left\{(f,g) \in L^p(X) \times L^p(X)\,:\,\exists r\in E(f),\,s\in
E(g);\,
\frac{1}{r}+\frac{1}{s} = \frac{1}{p}\,\right\}.$$
Evidently, $\Gamma_2 \subseteq \Gamma_1.$
If $ \mu(X)<\infty$ then $(L^p(X), \Gamma_2)$ is a partial *-algebra. But for
$\mu(X)=\infty$ the distributive property may fail (this is due to the fact
that  in the
case $\mu(X)=\infty$, the family of $L^p$-spaces form a lattice which do not
reduce to a chain). The set $\Gamma_2$ will be no longer considered here.

As shown before, for $p\geq 2,\; (L^p(X),\,C(X))$ is *-semisimple; then we
can define in $L^p(X)$ a {\em weak}-multiplication in the following way.
If $f,g \in L^p(X)$, we say that $f$ is a weak-multiplier of $g$, if there
exists a unique
element $h\in L^p(X)$ such that
$$\Omega(g\phi,f^*\psi)= \Omega(h\phi,\psi) \hspace{3mm} \forall \Omega \in
\Slp,\;\forall \phi,\psi  \in C(X).$$
It is worth remarking that, in this case, the uniqueness of $h$ follows from
Proposition \ref{rep} and so we do not need to require it.

\noindent Let $\Gamma_w$ denote the set of pairs $(f,g) \in L^p(X)\times
L^p(X)$ such that $f$ is a weak-multiplier of $g$. It is very easy to see that
$\Gamma_w=\Gamma_1$.

Another way to introduce a partial multiplication is to consider the so-called
{\em closable} elements.
In the following discussion we will not suppose that $X$ is compact (and so
$\mu (X)$ is not necessarily finite).
\begin{defn}
Let $f\in L^p(X)$ we say that $f$ is closable if the linear map
$$ T_f: \phi \in \ C_0(X) \mapsto f\phi \in L^p(X) $$ is closable as a densely
defined linear map in $L^p(X)$
\end{defn}
As is known \cite[\S 36]{Kot}, $T_f$ is closable if, and only if, it has an
adjoint $T_f^\prime$ whose domain $D(T_f^\prime)$ is weakly dense in
$L^{p'}(X)$.
It is easy to see that
\begin{equation}
 D(T_f^\prime)=\left\{ g \in L^{p'}(X) \,|\, f^*g \in L^{p'}(X) \right\}.
\label{dom}
\end{equation}
For $p>1$, this set is clearly weakly dense in $L^{p'}$, since it contains
$C_0(X)$.

\noindent Thus we have partially proved
\begin{prop}Let $1< p<\infty$. Every $f\in L^p(X)$ is closable.
If $\mu (X)<\infty$ then the statement holds also for $p=1$.
\end{prop}
\begin{pf}
It remains to prove the case $p=1$, under the assumption  $\mu (X)<\infty$.
This can be shown by an argument similar to that used in the first part of the
proof of Lemma \ref{PRD}
\end{pf}
For each element $g$ in the domain $D(\overline{T}_f)$ of the closure
$\overline{T}_f$ we can then
define as in \cite{bt2} a {\em strong} multiplication $\bullet$ setting
$$ f\bullet g  = \overline{T}_f g.$$
But, taking into account that $ \overline{T}_f$ coincides with the
double-adjoint of $T_f$,in analogy to \eqref{dom}, one has, for $p>1$
\begin{equation}
 D(\overline{T}_f)=\left\{ g \in L^{p}(X) \,|\, fg \in L^{p}(X) \right\}.
\label{dom2}
\end{equation}

If we denote with $\Gamma_s$ the set of pairs $(f,g) \in L^p(X)\times L^p(X)$
for which $f\bullet g$ is well-defined, from the above discussion it follows
that, for $p\geq 2$, $\Gamma_s=\Gamma_w$.
An interesting consequence of this fact is the following
\begin{cor}
Let $f,g \in L^p(X)$, $1< p<\infty$. The product $fg$ is in $L^p(X)$ if, and
only if, there
exists a sequence $g_n \in C_0(X)$ such that $g_n \stackrel{p}\longrightarrow
g$
and the sequence $fg_n$ converges in $L^p(X)$.
\end{cor}

We end this paper by remarking that the problem of refining the multiplication
in a quasi *-
algebra arise naturally from very simple examples (think of two step functions
in $(L^p(X,\mu),C_0(X))$). However, this problem has not always a positive
answer.
The relevant point here is that this is always possible, in non trivial way,
for *-semisimple CQ*-algebras \cite{bt2}.

\end{document}